\documentclass[11pt,a4paper]{amsart}
\usepackage{amsmath,amsfonts,amssymb,epsf,mathrsfs,xypic}
\usepackage[T1]{fontenc}
\usepackage[backref=page,bookmarksopen=true]{hyperref}
\usepackage{enumerate}
\numberwithin{equation}{section}

\let\bbbibitem\bibitem
\renewcommand{\bibitem}[2][]{\bbbibitem[#1]{#2}\label{#2}}

\def\fin{\hfill\hbox{\hskip .2cm $\square$}\medskip}

\newtheorem{theo}{Theorem}[section]
\newtheorem{lemma}[theo]{Lemma}
\newtheorem{prop}[theo]{Proposition}
\newtheorem{coro}[theo]{Corollary}

\newenvironment{demo}[1][\hspace{-3pt}]{{\noindent\em Proof #1.~ }}{\fin}

\def\cal{\mathcal}
\def\a{\alpha}
\def\b{\beta}

\def\g{\gamma}
\def\G{\Gamma}
\def\R{{\mathbb R}}

\def\Ker{{\rm Ker\,}}
\def\Im{{\rm Im\,}}
\def\C{{\mathbb C}}

\def\R{{\mathbb R}}
\def\fd{\longrightarrow}

\def\la{\langle}
\def\ra{\rangle}

\def\om{\omega}

\def\e{\varepsilon}

\def\L{{\cal L}}

\def\s{\sigma}
\def\l{\lambda}
\def\Th{\Theta}
\def\dt{\delta}
\def\rmo{\sqrt{-1}}
\def\dr{d'}
\def\db{d''}

\def\L{\Lambda}
\def\ext{\widetilde f_\infty}
\def\k{\kappa}

\def\E{{\cal E}}

\begin{document}

\title[Extensions with estimates of cohomology classes]{Extensions with estimates of cohomology classes}

\author{Vincent Koziarz} 
\address{IECN, Nancy-Universit\'e, CNRS, INRIA, Boulevard des Aiguillettes B. P. 70239, F-54506 
Vand\oe uvre-l\`es-Nancy, France}
\email{koziarz@iecn.u-nancy.fr}

\date{\today}

\sloppy

\begin{abstract}  We prove an extension theorem of ``Ohsawa-Takegoshi type'' for Dolbeault $q$-classes of cohomology ($q\geq 1$) on smooth compact hypersurfaces in a weakly pseudoconvex K\"ahler manifold.
\end{abstract}

\maketitle

\section{Introduction}

Let $Y$ be a complex submanifold of a K\"ahler manifold $X$ and let $L'$ be a Hermitian line bundle on $X$. First consider the following 

\medskip

\noindent {\bf Problem.} {\it Let $f$ be a smooth $D''$-closed section of $\Lambda^{0,q}T^\star_X\otimes L'$ over $Y$ satisfying a suitable $L^2$ condition. Can we find a smooth $D''$-closed extension $F$ of $f$ to $X$ together with a good $L^2$ estimate for $F$ on $X$?}

\medskip

The first result of this kind was obtained by T.~Ohasawa and K.~Takegoshi~\cite{OT} in the case when $Y$ is a hyperplane of a bounded pseudoconvex domain ~$X$ in $\C^n$, $L'$ is the trivial bundle and $q=0$. It was further generalized by L.~Manivel~\cite{Ma} (with a simplified proof by J.-P.~Demailly~\cite{De3}) in the following setting: $X$ is a weakly pseudoconvex manifold, $Y$ is the zero set of a holomorphic section of a rank $r$ Hermitian bundle over $X$, $L'=K_X\otimes L$ where $L$ is a Hermitian line bundle whose curvature satisfies appropriate positivity properties, $K_X$ is the canonical bundle of $X$, and $q=0$. When $q\geq 1$, the method leads to a new technical difficulty occurring in the regularity argument for $(0,q)$ forms. In \cite{De3}, Demailly suggests an approach to overcome this difficulty but, to our knowledge, the complete arguments did not appear anywhere. In this paper, we rather consider the

\medskip

\noindent {\bf Modified problem.} {\it Let $q\geq 1$ and $f$ be a smooth $D''$-closed section of $\Lambda^{0,q}T^\star_X\otimes L'$ over $Y$ satisfying a suitable $L^2$ condition. Can we find a smooth $D''$-closed extension $F$ of $f$ to $X$ as a cohomology class (i.e.$[F_{|Y}]=[f]\in H^q(Y,L')$) together with a good $L^2$ estimate for $F$ on $X$?}

\medskip

Observe that if $Y$ is a Stein submanifold (this happens e.g. when $X$ is a Stein manifold), the modified problem is not relevant when $q\geq 1$ since the Dolbeault group $H^q(Y,L')$ vanishes. In contrast, we will focus here on the case when $Y$ is a smooth compact hypersurface of a weakly pseudoconvex K\"ahler manifold $X$. This situation naturally happens, for example, when $X$ is a compact K\"ahler manifold, or when $X$ is a holomorphic family of projective algebraic manifolds fibered over the unit disc.

\begin{theo}\label{theop}
Let $(X,\om)$ be a weakly pseudonconvex $n$-dimensional K\"ahler manifold, and let $Y\subset X$ be the zero set of a holomorphic section $s\in H^0(X,E)$ of a Hermitian line bundle $(E,h_E)$; the subvariety $Y$ is assumed to be compact and nonsingular.
Let $L$ be a line bundle endowed with a smooth Hermitian metric $h_L$ such that
\begin{eqnarray}
&&\rmo\Theta(L)+\sqrt{-1}\dr\db \log|s|^2\geq 0,\label{assone}\\ 
&&\rmo\Theta(L)+\sqrt{-1}\dr\db\log|s|^2\geq \alpha^{-1}\rmo\Theta(E) \hbox{ for some }\alpha\geq 1,\label{asstwo}\\
&&|s|^2\leq e^{-\alpha}\label{assthree}
\end{eqnarray}
on $X$. Let $0<\k\leq 1$ and let $\Omega\subset X$ be a relatively compact open subset containing $Y$. Then, for any $q\geq0$ and every smooth $D''$-closed $(0,q)$-form $f$ with values in $K_X\otimes L$ over $Y$, there exists a smooth extension $F$ of $f$ to $\Omega$ as a cohomology class (i.e. $[F_{|Y}]=[f]\in H^q(Y,K_X\otimes L)$) such that
$$\int_{\Omega}\frac{|F|^2}{|s|^{2(1-\k)}}dV_\om\leq \frac{C}{\k}\int_Y\frac{|f|^2}{|ds|^2}dV_{Y,\om}
$$
where $C$ is a numerical constant depending only on $\Omega$, $E$, $L$ and $q$.
\end{theo}

The norm of the forms with values in bundles will always be computed with respect to the one induced by $\om$, $h_E$ and $h_L$. Also, $\Theta(E)$ (resp. $\Theta(L)$) will always denote the curvature of the Hermitian line bundle $(E,h_E)$ (resp. $(L,h_L)$). When the metrics will be twisted by some positive functions, the weights will appear explicitely in the formulae.

\bigskip

Our proof follows many of the ideas outlined in~\cite{De3}. First, 
using the weight bumping technique (and the adapted Bochner-Kodaira-Nakano inequality) initiated by Ohsawa and Takegoshi, for all $\e>0$, we build extensions of $f$ of class $C^1$ whose $L^2$ norm is controlled, and which are ``approximately'' $D''$-closed (in the sense that the $L^2$ norm of their $D''$-derivative is bounded by a constant times $\e$). Philosophically, passing to the limit as $\e\rightarrow 0$ should provide the desired extension but the limiting elliptic differential system is singular along $Y$ and this forbids the direct use of elliptic regularity arguments. Then, at this point, our strategy differs from Demailly's. Instead, we construct ``approximate'' $q$-cocycles $\zeta_\e$  in \v Cech cohomology corresponding to the previous extensions via an effective Leray's isomorphism, in a similar fashion as Y.-T.~Siu in \cite{Si}. During the process, we solve local $D''$-equations by standard techniques of L.~Hörmander.
Then, we can take the limit as $\e\rightarrow 0$ and use the ellipticity of the Laplacian in bidigree $(0,0)$ to ensure the smoothness of the extending cocycle $\zeta$. Finally, reversing the process, we get a smooth extension $F$ of $f$ as a cohomology class. 
Notice that the constant $C$ in Theorem~\ref{theop} is mainly related to a finite covering of $\Omega\supset Y$ by Stein open subsets (which is used to apply Leray's isomorphism) and the norm of the derivatives of a partition of unity subordinate to this finite covering. This explains in part why we need $Y$ to be compact.

\bigskip

A consequence of Theorem~\ref{theop} is a qualitative surjectivity theorem for restriction morphisms in Dolbeault cohomology:
\begin{coro}
Let $X$, $Y$, $E$ and $L$ be as in Theorem~\ref{theop} i.e. satisfying $(\ref{assone})$, $(\ref{asstwo})$ and $(\ref{assthree})$. Then the restriction morphism
$$H^q(X,K_X\otimes L)\fd H^q(Y,(K_X\otimes L)_{|Y})
$$
is surjective for any $q\geq 0$.
\end{coro}

Applying Theorem~\ref{theop} to $E=\C$ and to any semi-positive line bundle $L$ (for instance $L=\C$), we also easily get the following corollary which contains a special case of the invariance of the Hodge numbers for a family of compact K\"ahler manifolds (a result due to K.~Kodaira and D.~Spencer):
\begin{coro}Let $\pi:{\mathfrak X}\rightarrow \Delta$ be a proper holomorphic submersion over the unit disc and $L$ a semi-positive line bundle on ${\mathfrak X}$. Assume that ${\mathfrak X}$ is a K\"ahler manifold of dimension $n+1$. Then, for any $q\geq0$, $h^{n,q}(X_t,L):={\rm dim}\, H^{n,q}(X_t,L)$ is independent of $t\in\Delta$ (where $X_t=\pi^{-1}(t)$).
\end{coro}

\noindent{\bf Acknowledgments.} {I would like to thank Mihai P\u aun for many valuable discussions. I would also like to thank Jean-Pierre Demailly for explaining me details on his article~\cite{De3}, as well as Benoît Claudon and Dror Varolin for useful comments on an earlier version of this paper.}

\section{Preliminary material}
From now on, we assume that $X$, $Y$, $L$ and $E$ satisfy the hypotheses of Theorem~\ref{theop}.

Let $c\in\R$ such that $\overline{\Omega}\subset X_c:=\{x\in X\,,\,\psi(x)<c\}$, where $\psi$ is the plurisubharmonic exhaustion of $X$. Let ${\cal U}=\{U_j\}_{j\in J}$ be a finite covering of the closure of ${\Omega}$ by 
coordinate charts $\phi_j:B\fd U_j$ where $B$ is the unit ball in $\C^n$, and such that $U_j\subset {X_c}$ for all $j$. 
Denoting by $\nu$ the standard Hermitian norm on $\C^n$, we assume that the functions $\varphi_j:=\nu\circ\phi^{-1}_j$ satisfy
\begin{equation}\label{positiv}
\rmo\Theta(L)-\b\rmo\Theta(E)+\rmo\dr\db\varphi_j\geq \omega
\end{equation}
on $U_j$ for any $\b\in[0,1]$ (this is always possible if the $U_j$'s are chosen small enough). For any multi-index $(j_0,\dots,j_\ell)$, we shall denote by $\varphi_{j_0,\dots, j_\ell}$ the function $\sum_{i=0}^\ell\varphi_{j_i}$ which is defined on the intersection $U_{j_0,\dots,j_\ell}:=U_{j_0}\cap\dots\cap U_{j_\ell}$. 

If $F$ is a smooth Hermitian vector bundle over $X$ and $U\subset X$ is an open subset then, for all integer $k$, we denote by $\E^k(U,F)$ the space of sections of $F$ over $U$ which are of class $C^k$ and by $\E^k_c(U,F)$ those with compact support. We also denote by $W^k(U,F)$ the Sobolev space of sections whose derivatives (in the sense of distribution theory) up to order $k$ are in $L^2$.

Let us recall three useful results taken from~\cite{De1} (Remark 1.6, Lemma 3.3
and Lemma 6.9):
\begin{prop}\label{usef}
\begin{enumerate}[{\rm (}a{\rm )}]
\item $X_c\backslash Y$ is complete K\"ahler. 
\item Let $\om$ and $\om'$ be two Hermitian forms on $T_X$ such that $\om\leq\om'$. Let $E$ be a Hermitian vector bundle on $X$. Then, for any $q\geq 0$ and any $u\in\L^{n,q} T^\star_X\otimes E$, $|u|_{\om'}^2\, dV_{\om'}\leq |u|^2_\om\, dV_\om$.
\item Let $\Omega$ be an open subset of $\C^n$ and $Y$ a complex analytic subset of $\Omega$. Assume that $v$ is a $(p,q-1)$-form with $L^2_{\rm loc}$ coefficients and $w$ a $(p,q)$-form with $L^1_{\rm loc}$ coefficients such that $\db v=w$ on $\Omega\backslash Y$ (in the sense of distribution theory). Then $\db v=w$ on $\Omega$.
\end{enumerate}
\end{prop}

The following lemma is a consequence of a classical result (see \cite{De2}, Corollary 5.3):
\begin{lemma}\label{loc}
Let $m$ and $p$ be positive integers.
\begin{enumerate}[{\rm (}a{\rm )}]
\item Let $v\in W^m(U_{j_0,\dots,j_\ell},\Lambda^{n,p}T^\star_X\otimes L\otimes E^{-1})$ such that  $D''v=0$ and
$$\int_{U_{j_0,\dots,j_\ell}} |v|^2 e^{-\varphi_{j_0,\dots,j_\ell}}\, dV_{\omega}<+\infty.
$$
Then there exists a $(n,p-1)$ form $u\in W^{m+1}(U_{j_0,\dots,j_\ell},\Lambda^{n,p-1}T^\star_X\otimes L\otimes E^{-1})$ such that $D''u=v$ and
$$\int_{U_{j_0,\dots,j_\ell}} |u|^2 e^{-\varphi_{j_0,\dots,j_\ell}}\, dV_{\omega}\leq \frac{1}{p}\int_{U_{j_0,\dots,j_\ell}} |v|^2e^{-\varphi_{j_0,\dots,j_\ell}}\, dV_{\omega}.
$$
\item Let $0<\k\leq 1$ and $\e>0$. Let $v\in W^m(U_{j_0,\dots,j_\ell},\Lambda^{n,p}T^\star_X\otimes L)$ such that  $D''v=0$ and
$$\int_{U_{j_0,\dots,j_\ell}} \frac{|v|^2}{(|s|^2+\e^2)^{1-\k}}e^{-\varphi_{j_0,\dots,j_\ell}}\, dV_{\omega}<+\infty.
$$
Then there exists a $(n,p-1)$ form $u\in W^{m+1}(U_{j_0,\dots,j_\ell},\Lambda^{n,p-1}T^\star_X\otimes L)$ such that $D''u=v$ and
$$\int_{U_{j_0,\dots,j_\ell}} \frac{|u|^2}{(|s|^2+\e^2)^{1-\k}}e^{-\varphi_{j_0,\dots,j_\ell}}\, dV_{\omega}\leq \frac{1}{p}\int_{U_{j_0,\dots,j_\ell}} \frac{|v|^2}{(|s|^2+\e^2)^{1-\k}}e^{-\varphi_{j_0,\dots,j_\ell}}\, dV_{\omega}.
$$

\end{enumerate}
\end{lemma}

\begin{demo}
We only check the hypotheses of Corollary~5.3 in~\cite{De2}. The open subset $U_{j_0,\dots,j_\ell}\subset X$ is Stein and on $U_{j_0,\dots,j_\ell}$, the line bundle $L\otimes E^{-1}$, resp. $L$, endowed with its metric twisted by $e^{-\varphi_{j_0,\dots,j_\ell}}$, resp. $(|s|^2+\e^2)^{-(1-\k)}e^{-\varphi_{j_0,\dots,j_\ell}}$, has curvature

$$\rmo\Theta(L)-\rmo\Theta(E)+\rmo d'd''\varphi_{j_0,\dots,j_\ell},
$$
resp. 
$$\rmo\Theta(L)+(1-\k)\rmo d'd''\log(|s|^2+\e^2)+\rmo d'd''\varphi_{j_0,\dots,j_\ell},
$$
which, by inequality (\ref{bel}) and assumption~(\ref{positiv}), is bounded from below by
$$\rmo\Theta(L)-\rmo\Theta(E)+\rmo d'd''\varphi_{j_0}\geq \om,
$$
resp.
$$\rmo\Theta(L)-(1-\k)\frac{\la\rmo\Theta(E)s,s\ra}{|s|^2+\e^2}+\rmo d'd'' \varphi_{j_0}\geq \om.
$$
The fact that $u$ can be chosen in the Sobolev space $W^{m+1}$ comes from the ellipticity of the Laplacian. Let us explain why in the case $(a)$, the case $(b)$ being completely similar. In fact, the $D''$-equation is solved using complete metrics $\om_\e$ on $U_{j_0,\dots,j_\ell}$, such that $\om_\e\geq\om$ and $\om_\e\rightarrow\om$ as $\e\rightarrow 0$. For any $\e>0$, the corresponding minimal solution  $u_\e$ (i.e. the one satisfying $D'' u_\e=v$ and $u_\e\in (\Ker D'')^{\perp_{\om_\e}}$) is such that
$$\int_{U_{j_0,\dots,j_\ell}} |u_\e|_{\om_\e}^2 e^{-\varphi_{j_0,\dots,j_\ell}}\, dV_{\omega_\e}\leq \frac{1}{p}\int_{U_{j_0,\dots,j_\ell}} |v|_{\om_\e}^2e^{-\varphi_{j_0,\dots,j_\ell}}\, dV_{\omega_\e}\leq \frac{1}{p}\int_{U_{j_0,\dots,j_\ell}} |v|^2e^{-\varphi_{j_0,\dots,j_\ell}}\, dV_{\omega}$$
where the latter inequality comes from Proposition \ref{usef} $(b)$. Then, there exists a sequence $(\e_\mu)$ converging to 0 such that $u_{\e_\mu}$ converges weakly to some $u$ in $L^2_{\rm loc}$ as $\mu\rightarrow+\infty$: $u$ satisfies $D''u=v$,
$$\int_{U_{j_0,\dots,j_\ell}} |u|^2 e^{-\varphi_{j_0,\dots,j_\ell}}\, dV_{\omega}\leq \frac{1}{p}\int_{U_{j_0,\dots,j_\ell}} |v|^2e^{-\varphi_{j_0,\dots,j_\ell}}\, dV_{\omega}
$$
but also $u\in(\Ker D'')^{\perp_\om}=\overline{\Im (D'')^{\star_\om}}$ (and therefore $(D'')^{\star_\om} u=0$). Indeed, $L^2(\om)\subset L^2({\om_\e})$ because $|\,.\,|^2_{\om_\e}\, dV_{\om_\e}\leq |\,.\,|^2\, dV_{\om}$ and since $\om_\e\rightarrow\om$, we get $u\in(\Ker D'')^{\perp_\om}$ by the dominated convergence theorem. Finally, $u$ satisfies $D''u=v$ and $(D'')^\star u=0$ which is an elliptic differential system and standard arguments give $u\in W^{m+1}$ if $v\in W^m$.

Notice that we can skip the extraction of a weak limit if we only need a solution with the same estimate on a slightly smaller relatively compact open subset of $U_{j_0,\dots,j_\ell}$: all we have to do is take a complete metric on $U_{j_0,\dots,j_\ell}$ which coincides with $\om$ on the smaller subset.

\end{demo}

Finally, we also select  a smooth partition of unity $\{\sigma_j\}_{j\in J}$ subordinate to ${\cal U}$ (i.e. for each $j$, $\sigma_j\in\E^\infty_c(U_j,[0,1])$ and $\sum_{j\in J}\sigma_j(x)=1$ for any $x\in {\Omega}$).

\section{Proof of the theorem}

Recall that, by assumption, $Y\subset X$ is a smooth hypersurface so that $E\simeq{\cal O}_X(Y)$.

\medskip

\subsection{Construction of smooth extensions} In this section, we prove the following

\begin{lemma}\label{sext}
For any $k\geq 0$, there exists a smooth section
$$\widetilde f_\infty\in \E^\infty(X,\L^{n,q}T^\star_ X\otimes L)
$$
such that
\begin{enumerate}[{\rm (}a{\rm )}]
\item $\widetilde f_\infty$ coincides with $f$ in restriction to $Y$,
\item $|\widetilde f_\infty|=|f|$ at every point of $Y$,
\item $D''\widetilde f_\infty=0$ at every point of $Y$,
\item $s^{-1} D''\widetilde f_\infty\in \E^k(X,\L^{n,q+1}T^\star_ X \otimes L\otimes {\cal O}_X(-Y))$.
\end{enumerate}
\end{lemma}

\begin{demo}
Let us cover $Y$ by coordinate patches $W_j\subset X$ biholomorphic to polydiscs and with the following property:  if we denote the corresponding coordinates by $(z_j,w_j)\in\Delta\times\Delta^{n-1}$, where $w_j=(w_j^1,\dots, w_j^{n-1})$, then $W_j\cap Y=\{z_j=0\}$. On each $W_j$, we fix some holomorphic $\sigma_j\in \G(W_j,K_X\otimes L)$ which trivializes $K_X\otimes L$.

As explained in~\cite{De3}, the restriction map $(\L^{0,q}T^\star_X)_{|Y}\fd\L^{0,q}T^\star_Y$ can be viewed as an orthogonal projection onto a $C^\infty$ subbundle of $(\L^{0,q}T^\star_X)_{|Y}$. One might extend this subbundle from $W_j\cap Y$ to $W_j$ and then extend $f$ on ${W_j}$ by some smooth form $\widehat f_j\in \E^\infty(W_j,\L^{n,q}T^\star_ X\otimes L)$. Using a smooth partition of unity $\theta_j\in \E^\infty_c(W_j,\R)$, $\sum_j \theta_j=1$ on a neighbourhood of $Y$, we get a global smooth extension $\widehat f=\sum_j \theta_j\widehat f_j$ of $f$ which fulfills conditions $(a)$ and $(b)$. Since
$$(D''\widehat f)_{|Y}=(D''\widehat f_{|Y})=D'' f=0,$$
we can write $D''\widehat f=d\bar z_j\wedge g_j$ on $W_j\cap Y$ for some smooth $(0,q)$-forms $g_j$ which we extend arbitrarily to $W_j$. Then
$$\widetilde f_\infty:=\widehat f-\sum_j\theta_j\bar z_j g_j$$
coincides with $\widehat f$ on $Y$ and satisfies $(c)$.

We proceed by induction to get $(d)$. Assume that on each $W_j$,
$$D''\ext=z_j f_j(z_j,w_j) +\bar z_j^k \Bigl[d\bar z_j\wedge\sum_{|I|=q}a_I(w_j)\sigma_j d\bar w_j^I+\sum_{|I'|=q+1}b_{I'}(w_j)\sigma_j d\bar w_j^{I'}\Bigr]+\bar z_j^{k+1} h_j(z_j,w_j)
$$
for some $f_j,h_j\in \E^\infty(W_j,\L^{n,q+1}T^\star_ X\otimes L)$, $a_I,b_{I'} \in \E^\infty(\Delta^{n-1},\C)$, $k\geq 1$, and where the multi-indices $I$, $I'$ are increasing. We say that $\ext$ enjoys property $(P_k)$. Remark that such an equality implies that $s^{-1} D''\widetilde f_\infty\in \E^{k-2}(X,\L^{n,q+1}T^\star_ X\otimes L\otimes {\cal O}_X(-Y))$ if $k\geq 2$. Moreover, the extension $\widetilde f_\infty$ we just constructed satisfies property $(P_1)$ because $D''\widetilde f_\infty=0$ along $Y$.

Of course, $D''(D''\ext)=0$, but also the direct computation gives
$$D''(D''\ext)=z_jD''f_j(z_j,w_j)+k\bar z_j^{k-1} d\bar z_j\wedge\Bigl[\sum_{|I'|=q+1}b_{I'}(w_j)\sigma_j d\bar w_j^{I'}\Bigr]+\bar z_j^k h'_j(z_j,w_j)
$$
for some $h'_j\in\E^\infty(W_j,\L^{n,q+2}T^\star_ X\otimes L)$, hence 
the $b_{I'}$'s must vanish identically. So if we take
$$\widetilde f'_\infty=\ext-\sum_j\theta_j\frac{\bar z^{k+1}_j}{k+1}\sum_{|I|=q}a_I(w_j)\sigma_j d\bar w_j^I,
$$
we have
$$\begin{array}{rcl}
D''\widetilde f'_\infty &= & \displaystyle D''\ext-\sum_j \biggl(\frac{\bar z^{k+1}_j}{k+1}d''\theta_j+\theta_j \bar z^k_jd\bar z_j\biggr)\wedge\sum_{|I|=q}a_I(w_j)\sigma_j d\bar w_j^I\\
& &\hskip 4cm -\displaystyle\sum_j\theta_j\frac{\bar z^{k+1}_j}{k+1}D''\biggl(\sum_{|I|=q}a_I(w_j)\sigma_j d\bar w_j^I\biggr)\\
&=& \displaystyle \sum_j\theta_j\biggl(D''\ext-\bar z_j^k d\bar z_j\wedge\sum_{|I|=q}a_I(w_j)\sigma_j d\bar w_j^I\biggr)+\sum_j\bar z_j^{k+1} h''_j(z_j,w_j)\\
& = & \displaystyle\sum_j  z_j\theta_j f_j(z_j,w_j)+\sum_j \bar z_j^{k+1}\bigl(\theta_j h_j(z_j,w_j)+ h''_j(z_j,w_j)\bigr)
\end{array}
$$
for some $h''_j\in \E_c^\infty(W_j,\L^{n,q+1}T^\star_ X\otimes L)$. Then, $\widetilde f'_\infty$ enjoys property $(P_{k+1})$.

\end{demo}

\medskip

\subsection{Construction of approximate extensions with control}

Let $\theta:\R\fd [0,1]$ be a smooth function with support in $(-\infty,1)$, such that $\theta\equiv 1$ on $(-\infty,1/2]$ and $|\theta'|\leq 4$, and consider the truncated extension of $f$
$$\widetilde f_\e:=\theta(\e^{-2}|s|^2)\widetilde f_\infty$$
where $\ext$ is the extension provided by Lemma~\ref{sext}, such that $s^{-1} D''\widetilde f_\infty\in \E^k(X,\L^{n,q+1}T^\star_ X \otimes L\otimes {\cal O}_X(-Y))$ for some $k\geq 1$ which will be determined later. We wish to solve on $X$ the equation
$$D'' u_{\e}=D'' \widetilde f_{\e}$$
with estimate, and the additional constraint that  $u_{\e}$ vanishes along $Y$. We also expect some regularity on $u_\e$ in order to justify that $\widetilde f_{\e}-u_\e$ is a ($D''$-closed) extension of $f$. In general, we are not able to get this by the method we use here, and we can only produce approximate solutions.

\medskip

The fundamental tool is the following existence result (see \cite{De3}, \cite{Pa}):
\begin{theo}\label{ape}
Let $X$ be a complete K\"ahler manifold of dimension $n$ equipped with a (not necessarily complete) K\"ahler metric $\omega$, and let $L$ be a line bundle endowed with a smooth Hermitian metric. Assume that there exist two smooth bounded functions $\eta,\lambda>0$ on $X$ satisfying
\begin{equation}\label{cond}
\eta\sqrt{-1}\Theta(L)-\sqrt{-1}\dr\db\eta-\rmo\frac{\dr\eta\wedge\db\eta}{\lambda}\geq\sqrt{-1}\tau\dr\mu\wedge\db\mu
\end{equation}
for some positive function $\tau$ and some function $\mu$. Let us consider the (densely defined) modified $D''$ operators
$$Tu:=D''(\sqrt{\eta+\lambda}u)\ \ {\rm and}\ \ Su:=\sqrt\eta(D'' u)$$
acting on forms with values in $L$. Let $g=\db\mu\wedge g_0+g_2$ be a $L^2$ form of $(n,q+1)$ type ($q\geq 0$) with values in $L$ such that
\begin{enumerate}[{\rm (}a{\rm )}]
\item $D'' g=0$,
\item $g_0\in L^2 (X,\L^{n,q}T^\star_X\otimes L)$,
\item $C(g_0,\tau):=\int_X1/\tau|g_0|^2dV_\om<+\infty$,
\item $|g_2|^2\leq \gamma\, C(g_0,\tau)$ almost everywhere for some positive constant $\gamma$.
\end{enumerate}
Then, for any $u\in{\rm Dom}\,T^\star\cap{\rm Dom}\,S$, we have
$$\Bigl |\int_X \la g,u\ra dV_\om\Bigr |^2\leq C(g_0,\tau)\bigl(\|T^\star u\|^2+\|Su\|^2+\gamma\|u\|^2\bigr).$$
In particular, there exist $v\in L^2(X,\L^{n,q}T_X^\star\otimes L)$ and $w\in L^2(X,\L^{n,q+1}T_X^\star\otimes L)$ such that
$$Tv+\gamma^{1/2}w=g
$$
together with the estimate
$$\int_X |v|^2dV_\om+\int_X |w|^2dV_\omega\leq C(g_0,\tau).$$
\end{theo}

\bigskip

As before, let $c\in \R$ be such that $\overline{\Omega}\subset X_c$. For simplicity, we will assume in the sequel that $X=X_c$. We are going to apply Theorem~\ref{ape} to $D'' \widetilde f_{\e}$ on $X\backslash Y$. By Proposition~\ref{usef} $(a)$, $X\backslash Y$ can be equipped with a complete K\"ahler metric. As for the bundle $L$, we endow it with its original metric multiplied with the weight $|s|^{-2}$ in order to force the vanishing of the approximate solution along $Y$.

For any $\e>0$, we set $\sigma_\e:=-\log(\e^2+|s|^2)$. Remark that, because of the condition $|s|\leq e^{-\a}$, this function is positive for $\e$ small enough. Let $\chi:\R_+\fd\R_+$ be any strictly concave function whose derivative satisfies $1\leq\chi'\leq 2$ and such that $\chi(-\log(\e^2+e^{-2\a}))\geq 2\a$ for any $\e>0$ small enough (in \cite{De3} and \cite{Pa}, $\chi(t)=t+\log (1+t)$ but we will choose other functions picked in~\cite{MV}).  Let us define the two positive functions (again $\e$ is assumed to be small enough)
$$\eta_\e:=\chi(\sigma_\e)\ \ {\rm and}\ \ \lambda_\e:=-\frac{\chi'(\sigma_\e)^2}{\chi''(\sigma_\e)}.$$
Although this is done carefully in \cite{De3} and \cite{Pa}, we check quickly that $\eta_\e$ and $\lambda_\e$ fulfill condition $(\ref{cond})$ in Theorem~\ref{ape}. It is easy to see that
\begin{equation}\label{bel}
-\rmo\dr\db\s_\e\geq\rmo\frac{\e^2}{|s|^2}\dr\s_\e\wedge\db\s_\e-\frac{\la\rmo\Th(E)s,s\ra}{\e^2+|s|^2}
\end{equation}
and it is straightforward that
$$\dr\eta_\e=\chi'(\s_\e)\dr\s_\e\ ,\ \ \db\eta_\e=\chi'(\s_\e)\db\s_\e\ ,\ \ \dr\db\eta_\e=\chi'(\s_\e)\dr\db\s_\e+\chi''(\s_\e)\dr\s_\e\wedge\db\s_\e.$$
Thus, since $\chi'$ is positive,
$$-\rmo\dr\db\eta_\e\geq\Bigl(\frac{1}{\chi'(\s_\e)}\frac{\e^2}{|s|^2}+\frac{1}{\l_\e}\Bigr)\rmo\dr\eta_\e\wedge\db\eta_\e -\frac{\chi'(\s_\e)}{\e^2+|s|^2}\la\rmo\Th(E)s,s\ra.
$$
If $\e$ is small enough, then for any $x\in X$, $\eta_\e(x)\geq\chi(-\log(\e^2+e^{-2\a}))\geq 2\a$. Taking into account the curvature assumptions $(\ref{assone})$ and $(\ref{asstwo})$ in Theorem~\ref{theop} as well as the fact that $\chi'\leq 2$, we obtain
$$\eta_\e(\rmo\Th(L)+\rmo\dr\db\log|s|^2)\geq\frac{\eta_\e}{\a}\rmo\Th(E)\geq \frac{\chi'(\s_\e)}{\e^2+|s|^2}\la\rmo\Th(E)s,s\ra.
$$
Finally, summing up the two latter inequalities, we get
$$\eta_\e(\rmo\Th(L)+\rmo\dr\db\log|s|^2)-\rmo\dr\db\eta_\e-\frac{\rmo}{\l_\e}\dr\eta_\e\wedge\db\eta_\e\geq\frac{1}{\chi'(\s_\e)}\frac{\e^2}{|s|^2}\rmo\dr\eta_\e\wedge\db\eta_\e
$$
which proves that $(\ref{cond})$ is fulfilled with $\displaystyle\tau=\frac{1}{\chi'(\s_\e)}\frac{\e^2}{|s|^2}$ and $\mu=\eta_\e$. Now, we can write

$$D'' \widetilde f_\e=\db\eta_\e\wedge g_\e+\theta\Bigl(\frac{|s|^2}{\e^2}\Bigr)D''\ext
$$
where 
$$g_\e:=\Bigl(1+\frac{|s|^2}{\e^2}\Bigr)\theta'\Bigl(\frac{|s|^2}{\e^2}\Bigr)\frac{\widetilde f_\infty}{\chi'(\s_\e)}.
$$

A quick computation shows that 
\begin{equation}\label{lim}
\lim_{\e\rightarrow 0} \,\frac{1}{\e^2}\int_{X\backslash Y} \theta'\Bigl(\frac{|s|^2}{\e^2}\Bigr)^2 |\widetilde f_\infty|^2\,dV_{\omega}=c_0\int_{Y}\frac{|f|^2}{|ds|^2} dV_{Y,\om}
\end{equation}
for some ``universal'' constant $c_0$.
Therefore, since $\theta(\e^{-2}|s|^2)$ is supported in $\{|s|<\e\}$, and since $D''\ext=0$ on $Y$, for any $\gamma>0$,
$$\biggl|\theta\Bigl(\frac{|s|^2}{\e^2}\Bigr)D''\ext\biggr|\leq \gamma\int_{X\backslash Y} \frac{1}{\tau} \frac{|g_\e|^2}{|s|^2}dV_{\om}=\frac{\g}{\e^2}\int_{X\backslash Y} \Bigl(1+\frac{|s|^2}{\e^2}\Bigr)^2\theta'\Bigl(\frac{|s|^2}{\e^2}\Bigr)^2 |\widetilde f_\infty|^2\,dV_{\omega}
$$
if $\e>0$ is small enough.

Hence, we can apply Theorem~\ref{ape}: we find $u_{\e,\gamma}=\sqrt{\eta_\e+\l_\e} v_{\e,\gamma}$ and $w_{\e,\gamma}$ which satisfy the equation 
$$D'' u_{\e,\gamma}+\gamma^{1/2}w_{\e,\gamma}=D''\widetilde f_\e$$
on $X\backslash Y$ and such that

\begin{eqnarray}
\displaystyle\int_{X\backslash Y} \frac{|u_{\e,\gamma}|^2}{|s|^{2}(\eta_\e+\l_\e)}dV_{\omega} + \int_{X\backslash Y} \frac{|w_{\e,\gamma}|^2}{|s|^{2}}dV_{\omega}& \leq & \displaystyle\frac{1}{\e^2}\int_{X\backslash Y} \Bigl(1+\frac{|s|^2}{\e^2}\Bigr)^2\theta'\Bigl(\frac{|s|^2}{\e^2}\Bigr)^2\frac{|\widetilde f_\infty|^2}{\chi'(\s_\e)}dV_{\omega}\nonumber\\
& \leq &\displaystyle\frac{4}{\e^2}\int_{X\backslash Y} \theta'\Bigl(\frac{|s|^2}{\e^2}\Bigr)^2 |\widetilde f_\infty|^2\,dV_{\omega}.\label{esti}
\end{eqnarray}

\subsection{Regularization of the approximate solution}\label{regul}

Recall that $Y$ is a divisor such that $E\simeq {\cal O}_X(Y)$. Here we use a trick of Demailly:  we consider $s^{-1}u_{\e,\gamma}$ (resp. $s^{-1}w_{\e,\gamma}$) as a $L^2$ $(0,q)$-form (resp. $(0,q+1)$-form) with values in the twisted line bundle $K_X\otimes L\otimes{\cal O}_X(-Y)$ equipped with a smooth Hermitian metric. By Proposition~\ref{usef} $(c)$, we can write
$$D''(s^{-1}u_{\e,\gamma})+\gamma^{1/2}s^{-1}w_{\e,\gamma}=s^{-1}D''\widetilde f_\e
$$
not only on $X\backslash Y$ but also on $X$ because $s^{-1}u_{\e,\gamma}$ is locally $L^2$,  $s^{-1}w_{\e,\gamma}$ is $L^2$ hence locally $L^1$, and $s^{-1}D''\widetilde f_\e$ is of class $C^k$ hence locally $L^1$ (recall that $\ext$, as chosen in Lemma~\ref{sext}, is such that $s^{-1}D''\ext$ is of class $C^k$, $k\geq 1$). However, we do not know much about the regularity of $u_{\e,\g}$ and $w_{\e,\g}$.

But $\E^\infty_c(X,\Lambda^{n,q}T^\star_X\otimes L\otimes{\cal O}_X(-Y))$ is dense in ${\rm Dom}\,D''$ for the graph norm, where we consider $D''$ as an operator acting on $(n,q)$ forms on $X$ with values in $L\otimes{\cal O}_X(-Y)$. More precisely, the density holds when $X$ is endowed with a complete metric. If $X=X_c$ as we assumed above, we can work instead on $X_{c'}$ for some $c'>c$, and there exists on $X_{c'}$ some complete K\"ahler metric which coincides with $\om$ on $X_c$.

Then, we can find some $t_{\e,\gamma}\in \E^\infty(X,\L^{n,q}T^\star_{X}\otimes L\otimes{\cal O}_X(-Y))$, which is $L^2$, such that
\begin{equation}\label{control}
\left|\int_{X} \frac{|t_{\e,\gamma}|^2}{\eta_\e+\l_\e}dV_{\omega}-\int_{X} \frac{|s^{-1}u_{\e,\gamma}|^2}{\eta_\e+\l_\e}dV_{\omega}\right|\leq\e
\end{equation}
(recall that $\eta_\e$ is bounded by $2\alpha$ from below), and $D''(t_{\e,\gamma}-s^{-1}u_{\e,\gamma})$ has $L^2$ norm bounded by $\gamma^{1/2}$ from above. As a consequence,
$$D'' t_{\e,\gamma}=s^{-1}D''\widetilde f_\e+r_{\e,\gamma}
$$
on $X$, with $r_{\e,\gamma}\in \E^k(X,\L^{n,q+1}T^\star_X\otimes L\otimes E^{-1})$ since $D'' t_{\e,\gamma}$ and $s^{-1}D''\widetilde f_\e$ are of class $C^k$. Moreover, $r_{\e,\gamma}$ satisfies 
\begin{equation}\label{rgbound}
\int_X |r_{\e,\g}|^2\leq C_1^2\g
\end{equation}
for some positive constant $C_1$ depending on $f$, but not on $\e$ and $\g$ (see (\ref{esti}) and (\ref{lim})). Finally, $s^{-1}D''\widetilde f_\e$ is $D''$-closed on $X\backslash Y$, hence on $X$ by Proposition~\ref{usef} $(c)$, and therefore $D''r_{\e,\gamma}=0$.

\subsection{The choice of $\eta_\e$ and $\l_\e$}\label{choice}

Let us come now to the choice of $\eta_\e$ and $\l_\e$ (see \cite{MV} for more details). For any $0<\k\leq1$, we define for $t\geq 0$ the functions 
$$g_\k(t)=\k^{-1}e^{\k t}\ ,\ \  h_\k(t)=\int_0^t \frac{1}{2e^{\k y}-1} dy\ \ {\rm and}\ \ \chi_\k(t)=1+t+h_\k(t).$$
One checks immediatly that $1\leq\chi'_\k\leq 2$ and $\chi''_\k<0$. Moreover, 
$$\chi_\k(-\log(\e^2+e^{-2\a}))\geq 1- \log(\e^2+e^{-2\a})\geq 1+ 2\a -\log(1+\e^2 e^{2\a})\geq 2\a
$$
when $\e$ is small enough. Clearly, $\chi_\k(t)\leq 1+2t$ and it follows that
$$\frac{\chi_\k(t)}{g_\k(t)}\leq \frac{\k(1+2t)}{e^{\k t}}\leq 2 
$$
as is seen from a simple computation. Moreover,

$$-\frac{\chi'_\k(t)^2}{\chi''_\k(t)}\leq 2 g_\k(t).
$$
As a consequence, if we fix $\k\leq 1$ and take $\chi=\chi_\k$, we have
\begin{equation}\label{kbound}
\eta_\e+\l_\e\leq \frac{4}{\k(|s|^2+\e^2)^\k}.
\end{equation}

\subsection{Construction of $q$-cochains via Leray's isomorphism}

Recall that we have fixed a finite open covering ${\cal U}=\{U_j\}_{j\in J}$ of $\overline{\Omega}$. We endow the group $C^\ell_2({\cal U},\E^1(\Lambda^{n,p}T^\star_X\otimes L\otimes E^{-1}))$ of (alternate) $\ell$-cochains with values in $\E^1(\Lambda^{n,p}T^\star_X\otimes L\otimes E^{-1})$ which are $L^2$ with the norm

$$\|\varsigma^\ell\|^2=\max_{j_0<\dots<j_\ell}\int_{U_{j_0,\dots,j_{\ell}}} |\varsigma^\ell_{j_0,\dots,j_{\ell}}|^2 e^{-\varphi_{j_0,\dots,j_\ell}}\,dV_\om
$$
and for all $0<\k\leq 1$ and $\e>0$, we endow $C_2^\ell({\cal U},\E^1(\Lambda^{n,p}T^\star_X\otimes L))$ with the norm
$$\|\varsigma^\ell\|^2_{\k,\e}=\max_{j_0<\dots<j_\ell}\int_{U_{j_0,\dots,j_{\ell}}} \frac{|\varsigma^\ell_{j_0,\dots,j_{\ell}}|^2}{(|s|^2+\e^2)^{1-\k}} e^{-\varphi_{j_0,\dots,j_\ell}}\,dV_\om.
$$
Remark that in the case when $\varsigma^0$ is the $0$-cocycle associated to a section $\varsigma$ of $\E^1(\Lambda^{n,p}T^\star_X\otimes L\otimes E^{-1}))$ (resp. $\E^1(\Lambda^{n,p}T^\star_X\otimes L))$),
$$\|\varsigma^0\|^2\leq \int_{X} |\varsigma|^2 dV_\om\  \Bigl(\hbox{ resp. } \|\varsigma^0\|_{\k,\e}^2\leq \int_{X} \frac{|\varsigma|^2}{(|s|^2+\e^2)^{1-\k}}\, dV_\om\Bigr)
$$
since the $\varphi_j$'s are nonnegative.

Now, we construct a $(q+1)$-cocycle in $Z^{q+1}({\cal U},{\cal O}(K_X\otimes L\otimes E^{-1}))$ corresponding to $r_{\e,\gamma}$ via Leray's isomorphism between the Dolbeault and the \v Cech cohomology groups. In fact, we are mostly interested in the intermediate cochains which appear during the process and the control we have on their norm.
The extension $\ext$ of $f$ is supposed to be sufficiently regular (i.e. $k$ is large enough in Proposition~\ref{sext}) in order that $r_{\e,\g}$, and every cochain obtained by solving local $D''$-equations below, is at least of class $C^1$ (see section~\ref{regul}, Lemma~\ref{loc} and use Sobolev lemma: $W^m(U)\subset \E^1(U)$ for any open subset $U\subset X$ if $m>1+\frac{n}{2}$).

For notational simplicity, we denote $r^{\ell}_{\e,\gamma,(j_0,\dots,j_\ell)}\in \G(U_{j_0,\dots,j_\ell},\E^1(\Lambda^{n,q-\ell}T^\star_X\otimes L\otimes E^{-1}))$ by $r^{\ell}_{j_0,\dots,j_\ell}$.

First, we solve the equation $D'' r_j^0=r_{\e,\gamma}$ on the $U_j$'s, then we solve the equations
$$D'' r^{\ell+1}_{j_0,\dots,j_{\ell+1}}=(\delta r^\ell)_{j_0,\dots,j_{\ell+1}}\ \ {\rm on}\ \ U_{j_0}\cap\dots\cap U_{j_{\ell+1}}\ \ (0\leq\ell\leq q-1)
$$
using each time Lemma~\ref{loc} $(a)$. Finally, $\delta r^q_{\e,\gamma}\in Z^{q+1}({\cal U},{\cal O}(K_X\otimes L\otimes E^{-1}))$ is a representative in \v Cech cohomology of $[r_{\e,\g}]\in H^{q+1}(X,K_X\otimes L\otimes E^{-1})$.

\begin{lemma}\label{rbound}
For any $\ell$, we have
$$\|r^\ell_{\e,\g}\|^2\leq \frac{(\ell+1)\dots 2.1}{(q+1). q \dots(q-\ell+1)}\,\|r_{\e,\g}\|^2\leq \frac{(\ell+1)\dots 2.1}{(q+1). q \dots(q-\ell+1)}\,\int_{X} |r_{\e,\gamma}|^2 dV_\om.$$
\end{lemma}

\begin{demo}
For any $\ell\geq 1$,
$$\begin{array}{rcl}
\|r^\ell_{\e,\g}\|^2 &\leq &  \displaystyle  \frac{1}{q-\ell+1} \|\dt r^{\ell-1}_{\e,\g}\|^2\\
&\leq & \displaystyle  \frac{\ell+1}{q-\ell+1} \|r^{\ell-1}_{\e,\g}\|^2\\
\end{array}
$$
and
$$
\|r^0_{\e,\g}\|^2 \leq \frac{1}{q+1} \|r_{\e,\g}\|^2\leq \frac{1}{q+1}\int_{X} |r_{\e,\gamma}|^2 dV_\om
$$
by the estimate in Lemma~\ref{loc} $(a)$.

\end{demo}
\medskip

In a similar manner, we produce a $q$-cochain $\zeta_{\e,\g}\in C^{q}({\cal U},{\cal O}(K_X\otimes L))$ corresponding to the ``approximately'' $D''$-closed extension $\widetilde f_\e-st_{\e,\gamma}\in\G(X,\E^1(\Lambda^{n,q}T^\star_X\otimes L))$ of $f$. More precisely, on the $U_j$'s, we solve the equation $D''h^0_j=\widetilde f_\e-st_{\e,\gamma}+sr^0_{\e,\g,j}$, then we solve
$$D''h^{\ell+1}_{j_0,\dots,j_{\ell+1}}=(\delta h^\ell)_{j_0,\dots,j_{\ell+1}}+(-1)^{\ell+1} sr^{\ell+1}_{j_0,\dots,j_{\ell+1}}\ \ {\rm on}\ \ U_{j_0}\cap\dots\cap U_{j_\ell+1}\ \ (0\leq\ell\leq q-2)
$$
using Lemma~\ref{loc} $(b)$. This is indeed possible since the right-hand side is $D''$-closed: if
$$D''h^{\ell}_{j_0,\dots,j_{\ell}}=(\delta h^{\ell-1})_{j_0,\dots,j_{\ell}}+(-1)^{\ell} sr^{\ell}_{j_0,\dots,j_{\ell}}
$$
then
$$\begin{array}{rcl}
\displaystyle D''(\delta h^{\ell})_{j_0,\dots,j_{\ell+1}} =  (\delta(D'' h^\ell))_{j_0,\dots,j_{\ell+1}}& =& \displaystyle (\delta^2h^{\ell-1})_{j_0,\dots,j_{\ell+1}}+(-1)^\ell s(\delta r^\ell)_{j_0,\dots,j_{\ell+1}}\\
& = & \displaystyle 0+(-1)^\ell sD''r^{\ell+1}_{j_0,\dots,j_{\ell+1}}.
\end{array}
$$
Finally, let $\zeta_{\e,\g}:=\delta h_{\e,\g}^{q-1}+(-1)^{q}sr^q_{\e,\g}\in C^{q}({\cal U},{\cal O}(K_X\otimes L))$ (in particular $D''\zeta_{\e,\g}=0$ and $\zeta_{\e,\g}$ is actually smooth by ellipticity of $D''$ in bidegree $(0,0)$). 

\medskip

In the next proposition, we denote by ${\cal V}$ the finite Stein covering $\{V_j\}_{j\in J}=\{U_j\cap Y\}_{j\in J}$ of $Y$.

\begin{prop}\label{estim}
 Let $0<\k\leq 1$. The cochain $\zeta_{\e,\g}$ enjoys the following properties:
\begin{enumerate}[{\rm (}a{\rm )}]
\item 
$$\|\zeta_{\e,\g}\|^2_{\k,\e}\leq \frac{16(q+1)c_0}{\k}\int_Y\frac{|f|^2}{|ds|^2} dV_{Y,\om}+\beta(\e,\g)
$$
where $\b$ is a positive function such that $\b(\e,\g)\rightarrow 0$ as $\e,\g\rightarrow 0$ (recall that for any $\g>0$, $\zeta_{\e,\g}$ only exists if $\e>0$ is small enough).
\item For any $\e,\g$, ${\zeta_{\e,\g}}_{|Y}\in Z^q({\cal V},{\cal O}((K_X\otimes L)_{|Y}))$ is a representative in \v Cech cohomology of the cohomology class $[f]\in H^q(Y,(K_X\otimes L)_{|Y})$.

\end{enumerate}
\end{prop}

\begin{demo}
Let $\chi=\chi_\k$ be as in section~\ref{choice}. We use the corresponding $\eta_\e$ and $\l_\e$. 

\begin{enumerate}[{\rm (}\it a{\rm )}] 
\item We have 
$$\begin{array}{rcl}
\|\zeta_{\e,\g}\|_{\k,\e} & \leq & \|\dt h^{q-1}_{\e,\g}\|_{\k,\e}+\|sr^q_{\e,\g}\|_{\k,\e}\\
&\leq&\displaystyle \sqrt{q+1}\|h^{q-1}_{\e,\g}\|_{\k,\e}+e^{-\frac{\a\k}{2}}\|r_{\e,\g}\|
\end{array}
$$
by Lemma~\ref{rbound}, since $|s|^2(|s|^2+\e^2)^{-(1-\k)}\leq |s|^\k\leq e^{-\frac{\a\k}{2}}$ according to assumption $(\ref{assthree})$. In the same way, for any $\ell\geq 1$,
$$\begin{array}{rcl}
\|h^\ell_{\e,\g}\|_{\k,\e} & \leq & \displaystyle\frac{1}{\sqrt{q-\ell}}(\|\dt h^{\ell-1}_{\e,\g}\|_{\k,\e}+\|sr^\ell_{\e,\g}\|_{\k,\e})\\
&\leq&\displaystyle\frac{\sqrt{\ell+1}}{\sqrt{q-\ell}} \|h^{\ell-1}_{\e,\g}\|_{\k,\e}+ \frac{\sqrt{(\ell+1)\dots 2.1}}{\sqrt{(q+1).q\dots(q-\ell)}}e^{-\frac{\a\k}{2}}\|r_{\e,\g}\|
\end{array}
$$
where we also used the estimate in Lemma~\ref{loc} $(b)$. Finally,
$$\begin{array}{rcl}
\|h^0_{\e,\g}\|_{\k,\e} & \leq & \displaystyle\frac{1}{\sqrt{q}}(\|\widetilde f_\e-st_{\e,\gamma}\|_{\k,\e}+\|sr^0_{\e,\g}\|_{\k,\e})\\
&\leq&\displaystyle\frac{1}{\sqrt{q}}\|\widetilde f_\e-st_{\e,\gamma}\|_{\k,\e}+ \frac{1}{\sqrt{(q+1).q}}e^{-\frac{\a\k}{2}}\|r_{\e,\g}\|.
\end{array}
$$
Collecting all these inequalities, we get
$$\|\zeta_{\e,\g}\|_{\k,\e}\leq \sqrt{q+1}\|\widetilde f_\e-st_{\e,\gamma}\|_{\k,\e}+q e^{-\frac{\a\k}{2}}\|r_{\e,\g}\|\leq \sqrt{q+1}\|\widetilde f_\e-st_{\e,\gamma}\|_{\k,\e}+q\,C_1\g^{1/2}.
$$
Now, it is easy to see that
$$\int_{X} \frac{|\widetilde f_\e-st_{\e,\gamma}|^2}{(|s|^2+\e^2)^{1-\k}}\, dV_\om\leq \int_{X} \frac{|st_{\e,\gamma}|^2}{(|s|^2+\e^2)^{1-\k}}\, dV_\om+C_2\,\e^{2\k}\\
$$
for some constant $C_2$, as $\widetilde f_\e$ is uniformly bounded with support in $\{|s|<\e\}$. Finally, by $(\ref{kbound})$,
$$\int_{X} \frac{|st_{\e,\gamma}|^2}{(|s|^2+\e^2)^{1-\k}}\, dV_\om\leq\int_{X} (|s|^2+\e^2)^\k |t_{\e,\g}|^2\, dV_\om\leq \frac{4}{\k}\int_{X}\frac{|t_{\e,\g}|^2}{\eta_\e+\l_\e}\, dV_\om
$$
and the desired inequality follows from (\ref{esti}), (\ref{control}) and (\ref{lim}).

\item It is clear since on $U_{j_0,\dots,j_\ell+1}\cap Y$, the restriction of $(\delta h^\ell)_{j_0,\dots,j_\ell}$ is always $D''$-closed, hence we construct an ``exact'' representative in restriction to $Y$.
\end{enumerate}
\end{demo}

\subsection{Passing to the limit and reversing the process to get the extension}

Now, we just have to make $\g$ and $\e$ go to zero and extract a weak limit of $\zeta_{\e,\g}$. This weak limit $\zeta$ is an element of $Z^{q}({\cal U},{\cal O}(K_X\otimes L))$ since $D''\zeta_{\e,\g}=0$ for any $\e,\g$, and $\delta \zeta_{\e,\g}=(-1)^q\delta (sr^q_{\e,\gamma})$ which, by $(\ref{rgbound})$ and Proposition~\ref{rbound}, has $L^2$ norm bounded by some constant times $\g^{1/2}$, thus $\delta\zeta=0$. Moreover, $\zeta_{|Y}\in Z^q({\cal V},{\cal O}((K_X\otimes L)_{|Y}))$ is a representative in \v Cech cohomology of the cohomology class $[f]\in H^q(Y,(K_X\otimes L)_{|Y})$ since this is the case of any ${\zeta_{\e,\g}}_{|Y}$ by Proposition~\ref{estim} $(b)$. For all $j_0,\dots,j_q$, as the $\varphi_j$'s are bounded from above by 1, we obtain
$$\displaystyle\int_{U_{j_0,\dots,j_{q}}} \frac{|\zeta|^2}{|s|^{2(1-\k)}} \,dV_\om\leq \frac{16 e^{q+1}(q+1)c_0}{\k}\int_Y\frac{|f|^2}{|ds|^2} dV_{Y,\om}
$$
if we take the limit in the inequality of Proposition~\ref{estim} $(a)$.

Finally, we construct the desired extension $F$ in the following way.
For $0\leq \ell\leq q-1$, we produce $\xi^\ell\in C^{\ell}({\cal U},\E^\infty(\Lambda^{n,q-\ell-1}T^\star_X\otimes L))$ such that
\begin{enumerate}[$(i)$]
\item $(\delta\xi^{q-1})_{j_0,\dots, j_q}=\zeta_{j_0,\dots, j_q}$ on $U_{j_0}\cap\dots\cap U_{j_q}$,
\item $(\delta\xi^{\ell})_{j_0,\dots, j_{\ell+1}}=D''\xi^{\ell+1}_{j_0,\dots, j_{\ell+1}}$ on $U_{j_0}\cap\dots\cap U_{j_{\ell+1}}$ ($0\leq \ell\leq q-2$).
\end{enumerate}
These $\dt$-equations are solved by using the partition of unity $\{\sigma_j\}_{j\in J}$ subordinate to ${\cal U}$ in the following way:

$$\begin{array}{rcl}
\xi^{q-1}_{j_0,\dots, j_{q-1}} & = & \sum_i \sigma_i\, \zeta_{i,j_0,\dots, j_{q-1}}\,,\\
\xi^{\ell}_{j_0,\dots, j_{\ell}} & = & \sum_i \sigma_i\,D''\xi^{\ell+1}_{i,j_0,\dots, j_{\ell}}\ \ \ (0\leq \ell\leq q-2).
\end{array}$$
Finally, we set 
$$F=D''\xi^{0}_{j}=q!\sum_{j_1<\dots<j_{q} \atop j_i\not = j} \zeta_{j_{q},\dots, j_{1},j}\, d''\sigma_{j_1}\wedge\dots\wedge d''\sigma_{j_q}$$
on $U_{j}$. Then, $F$ defines a $D''$-closed section in $\G(\Omega,\E^\infty(\Lambda^{0,q}T^\star_X\otimes K_X\otimes L))$ such that $[F_{|Y}]=[f]\in H^q(Y,(K_X\otimes L)_{|Y})$ and the estimate
$$\int_{\Omega}\frac{|F|^2}{|s|^{2(1-\k)}} \,dV_\om\leq \frac{16 e^{q+1}(q+1)c_0 C_\sigma |J|!}{\k(|J|-q-1)!}\int_Y\frac{|f|^2}{|ds|^2} dV_{Y,\om}
$$
holds for some constant $C_\sigma$ depending only on the partition of unity $\{\sigma_j\}_{j\in J}$ and $q$ (one can take $C_\sigma=\max_{j\in J} |d''\sigma_j|^{2q}$).

\end{document}